\documentclass[11pt,a4paper]{article}
\usepackage[T1]{fontenc}
\usepackage[utf8]{inputenc}
\usepackage{lmodern}
\usepackage[a4paper,margin=27mm]{geometry}
\usepackage{amsmath,amssymb}
\usepackage{graphicx}
\usepackage{microtype}
\usepackage{caption}
\usepackage{needspace,placeins}
\usepackage{xurl}
\usepackage[hidelinks]{hyperref}
\hypersetup{pdftitle={Cancellation profiles of Fibonacci-Lucas zeta sums},pdfauthor={Payam Danesh}}
\title{Cancellation profiles of Fibonacci-Lucas zeta sums}
\author{Payam Danesh}
\date{}
\begin{document}
\maketitle
\begin{abstract}
Finite Fibonacci and Lucas sums involving zeta values at negative
integers can have individual terms far larger than their total. In this
paper, we study the distribution of these term magnitudes as the even
degree increases. For non-negative weights of at most exponential
growth, we prove a finite exponential expansion for the normalized
absolute-term sum with an explicit remainder at every fixed order. The
normalized magnitudes converge in total variation to a discrete
distribution determined by the even exponential generating function of
the weights. For the Fibonacci and Lucas sums, this distribution has its
unique maximum at index eight and we prove that the same index gives the
largest finite-sum term at every even degree at least twenty-six. We
also derive the factorial growth of the summation condition number. The
signed identities are placed within the established Bernoulli-polynomial
reflection framework and extended to a trace-one recurrence family,
including repeated roots and a vanishing scale. Exact rational
calculations and precision-refined evaluations illustrate the proved
bounds and distinguish a fixed cancellation profile from the growing
length of the sum.

\end{abstract}


\hypertarget{introduction-and-background}{%
\section*{1 Introduction and
background}\label{introduction-and-background}}

An exact identity determines the value of a sum without necessarily
describing the size or location of its individual terms. This
distinction is pronounced when Bernoulli numbers occur with binomial
coefficients and recurrence sequences. At moderate degree, the terms can
exceed their signed total by many orders of magnitude. The relation
between Bernoulli numbers and zeta values at negative integers belongs
to classical analytic number theory field {[}1{]}. Sun studied
combinatorial identities through dual sequences {[}2{]}. Frontczak
related Fibonacci and Bernoulli numbers through balancing polynomials
{[}3{]} and later in further studies, they gave Bernoulli-polynomial
convolution formulas containing the Fibonacci and Lucas identities
{[}4{]}, {[}5{]}. More recently, Adegoke treated the underlying
binomial-transform symmetries in a general sequence frame {[}6{]} . A
separate line of the work concerns the sensitivity and asymptotic
behavior of finite sums. Higham's analysis of floating-point summation
identifies the ratio of the absolute-term sum to the absolute signed sum
as the relevant component wise condition number {[}7{]}. Coefficient
asymptotic for meromorphic generating functions provide a general
framework for successive exponentially separated contributions {[}8{]}.
Building on the finite identities {[}9{]}, we use the positive Dirichlet
series for the zeta function to obtain such contributions directly, with
explicit bounds under a simple growth assumption on the weights.

The central objects here are the absolute-term sum and the distribution
obtained by dividing each term magnitude by that sum. For any
non-negative weight sequence with a positive second coefficient and at
most exponential growth, we prove a finite expansion whose coefficients
are values of the even exponential generating function. We then prove
convergence of the normalized magnitudes in total variation with an
explicit exponential rate. These statements do not require a vanishing
Bernoulli convolution, when applied to the Fibonacci and Lucas weights,
they identify a unique limiting largest index and yield a uniform
finite-degree conclusion. The largest term occurs at index eight for
every even degree at least twenty-six. Thus the strongest cancellation
remains near a fixed index while the number of summands increases. The
formula supply the algebraic cancellation mechanism; the results below
quantify the entire absolute-term profile and locate its largest member.

\hypertarget{the-algebraic-origin-of-the-finite-identities}{%
\section*{2 The algebraic origin of the finite
identities}\label{the-algebraic-origin-of-the-finite-identities}}

\hypertarget{reflection-and-the-coefficient-identity}{%
\subsection*{2.1 Reflection and the coefficient
identity}\label{reflection-and-the-coefficient-identity}}

All sequences in the algebraic part have complex entries and are indexed
by non-negative integers. We write
\(\mathbb{C\lbrack\lbrack}z\rbrack\rbrack\) for the ring of formal power
series in an indeterminate \(z\). Equality in this ring means equality
of every coefficient, without an assumption of analytic convergence. An
empty sum is zero and a zeroth power is one, including at a zero base
when it occurs in a polynomial expression.

The Bernoulli numbers \(B_{m}\) and Bernoulli polynomials \(B_{m}(t)\)
are defined by

\begin{equation}
\frac{z}{e^{z} - 1} = \sum_{m = 0}^{\infty}B_{m}\frac{z^{m}}{m!},\quad\quad\frac{ze^{tz}}{e^{z} - 1} = \sum_{m = 0}^{\infty}B_{m}(t)\frac{z^{m}}{m!}.
\tag{1}\label{eq:1}
\end{equation}

Here the quotients are interpreted by their power-series expansions at
zero. In particular, \(B_{m} = B_{m}(0)\), \(B_{0} = 1\),
\(B_{1} = - 1/2\), and \(B_{2r + 1} = 0\) for every integer
\(r \geq 1\). Analytically, the expansions in (1) converge for
\(|z| < 2\pi\). The standard finite expansion and reflection identities
are

\begin{equation}
B_{N}(t) = \sum_{j = 0}^{N}\binom{N}{j}B_{N - j}t^{j},\quad\quad B_{N}(1 - t) = ( - 1)^{N}B_{N}(t).
\tag{2}\label{eq:2}
\end{equation}

The function \(\zeta(s)\) denotes the meromorphic continuation of the
series \(\sum_{m \geq 1}^{}m^{- s}\), initially defined for \(Res > 1\).
With the Bernoulli convention in (1), the special values required here
are

\begin{equation}
\zeta(1 - m) = - \frac{B_{m}}{m}\quad(m \geq 2),\quad\quad\zeta(0) = - \frac{1}{2}.
\tag{3}\label{eq:3}
\end{equation}

The first formula in (3) deliberately excludes \(m = 1\). All zeta
arguments in the finite sums below are negative integers and no infinite
sum is assigned a regularized value. For a complex parameter \(h\) and a
sequence \(a = \left( a_{j} \right)_{j \geq 0}\), we define its
exponential generating function \(A\) and its scaled Bernoulli
convolution \(T_{N}^{(h)}\) by

\begin{equation}
A(z) = \sum_{j = 0}^{\infty}a_{j}\frac{z^{j}}{j!},\quad\quad T_{N}^{(h)}(a) = \sum_{j = 0}^{N}\binom{N}{j}h^{N - j}B_{N - j}a_{j}.
\tag{4}\label{eq:4}
\end{equation}

The scaling kernel is the formal series
\(K_{h}(z) = \sum_{m \geq 0}^{}B_{m}h^{m}z^{m}/m!\). It equals
\(hz/\left( e^{hz} - 1 \right)\) for \(h \neq 0\), while
\(K_{0}(z) = 1\). Its constant coefficient is one, so it is invertible
in \(\mathbb{C\lbrack\lbrack}z\rbrack\rbrack\). Multiplying the
generating functions in (4) gives us

\begin{equation}
\sum_{N = 0}^{\infty}T_{N}^{(h)}(a)\frac{z^{N}}{N!} = K_{h}(z)A(z),\quad\quad K_{h}( - z) = e^{hz}K_{h}(z).
\tag{5}\label{eq:5}
\end{equation}

The second relation in Eq.(5) is the source of the parity cancellation.

\textbf{Lemma 1.} For every \(h\mathbb{\in C}\), the following
conditions are equivalent: \(T_{2k}^{(h)}(a) = 0\) for all integers
\(k \geq 0\); \(e^{hz}A( - z) = - A(z)\) and \(e^{- hz/2}A(z)\) is an
odd formal power series. Each condition implies \(a_{0} = 0\).

\textbf{Proof.} A formal series is odd precisely when its values at
\(z\) and \(- z\) sum to zero. By (5), the even coefficients of
\(K_{h}A\) vanish exactly when

\begin{equation}
K_{h}(z)\lbrack A(z) + e^{hz}A( - z)\rbrack = 0.
\tag{6}\label{eq:6}
\end{equation}

Multiplication of (6) by the inverse of \(K_{h}\) proves the first
equivalence. By Putting \(H(z) = e^{- hz/2}A(z)\), the equation
\(H( - z) = - H(z)\) is equivalent to \(e^{hz}A( - z) = - A(z)\) after
multiplication by \(e^{hz/2}\). Finally, the constant term of that
equation gives \(2a_{0} = 0\).

At \(h = 1\), this is the generating-function form of the
anti-self-inverse condition for the alternating binomial transform
discussed in {[}2,6{]}.

To retain the explicit coefficient route, we define the Fibonacci
sequence \(\left( F_{j} \right)\) and Lucas sequence
\(\left( L_{j} \right)\) by

\begin{equation}
\begin{aligned}
 & F_{0} = 0,\quad F_{1} = 1,\quad F_{j + 2} = F_{j + 1} + F_{j}, \\
 & L_{0} = 2,\quad L_{1} = 1,\quad L_{j + 2} = L_{j + 1} + L_{j}.
\end{aligned}
\tag{7}\label{eq:7}
\end{equation}

The recurrence in (7) has characteristic roots and Binet formulas

\begin{equation}
\alpha = \frac{1 + \sqrt{5}}{2},\quad\beta = \frac{1 - \sqrt{5}}{2},\quad\quad F_{j} = \frac{\alpha^{j} - \beta^{j}}{\sqrt{5}},\quad L_{j} = \alpha^{j} + \beta^{j}.
\tag{8}\label{eq:8}
\end{equation}

The identities \(\alpha + \beta = 1\), \(\alpha^{2} = \alpha + 1\), and
\(\beta^{2} = \beta + 1\) follow from the characteristic equation.
Hence, for \(\mathcal{F}(z) = \sum_{j \geq 0}^{}F_{j}z^{j}/j!\),

\begin{equation}
e^{z}\mathcal{F}(z) = \sum_{j = 0}^{\infty}F_{2j}\frac{z^{j}}{j!},\quad\quad e^{z}\mathcal{F}( - z)\mathcal{= - F}(z).
\tag{9}\label{eq:9}
\end{equation}

In fact, the two exponents in \(e^{z}\mathcal{F}(z)\) are
\(\alpha + 1 = \alpha^{2}\) and \(\beta + 1 = \beta^{2}\), whereas those
in \(e^{z}\mathcal{F}( - z)\) are \(1 - \alpha = \beta\) and
\(1 - \beta = \alpha\). Subtracting the odd parts in (9) therefore gives
us

\begin{equation}
2e^{z}\sum_{m = 0}^{\infty}F_{2m + 1}\frac{z^{2m + 1}}{(2m + 1)!} = \sum_{j = 1}^{\infty}\left( F_{j} + F_{2j} \right)\frac{z^{j}}{j!}.
\tag{10}\label{eq:10}
\end{equation}

The remaining step is a finite coefficient calculation.

\textbf{Proposition 2.} For each integer \(N \geq 1\):

\begin{equation}
\begin{aligned}
 & \sum_{m = 0}^{\lfloor(N - 1)/2\rfloor}\binom{N}{2m + 1}F_{2m + 1}\left( 2 - 2^{N - 2m - 1} \right)B_{N - 2m - 1} \\
 & = \sum_{j = 1}^{N}\binom{N}{j}\left( F_{j} + F_{2j} \right)2^{N - j - 1}B_{N - j}.
\end{aligned}
\tag{11}\label{eq:11}
\end{equation}

\textbf{Proof.} The required kernel relation is

\begin{equation}
\frac{z}{e^{z} - 1} + \frac{z}{e^{z} + 1} = \frac{2ze^{z}}{e^{2z} - 1}.
\tag{12}\label{eq:12}
\end{equation}

By Multiplying (10) by \(z/\left( e^{2z} - 1 \right)\) and apply (12),
from (1) and the decomposition
\(z/\left( e^{z} + 1 \right) = z/\left( e^{z} - 1 \right) - 2z/\left( e^{2z} - 1 \right)\),
the two resulting kernels have expansions

\begin{equation}
\begin{aligned}
\frac{z}{e^{z} - 1} + \frac{z}{e^{z} + 1} & = \sum_{r = 0}^{\infty}\left( 2 - 2^{r} \right)B_{r}\frac{z^{r}}{r!}, \\
\frac{z}{e^{2z} - 1} & = \sum_{r = 0}^{\infty}2^{r - 1}B_{r}\frac{z^{r}}{r!}.
\end{aligned}
\tag{13}\label{eq:13}
\end{equation}

The coefficient of \(z^{N}/N!\) in a product of exponential generating
functions is the binomial convolution of their coefficients. On the left
it selects \(r = N - 2m - 1\) and on the right it selects \(r = N - j\).
Substitution in (13) gives us exactly (11) and every coefficient
involves a finite sum.

When \(N\) is even, all Bernoulli indices on the left of (11) are odd.
Indices greater than one contribute zero and the remaining index one has
multiplier \(2 - 2^{1} = 0\). Thus the right-hand side vanishes. The
next construction places this observation in a recurrence family and
separates its two endpoint terms.

\hypertarget{the-recurrence-family}{%
\subsection*{\texorpdfstring{2.2 The recurrence family
}{2.2 The recurrence family }}\label{the-recurrence-family}}

By letting \(q\mathbb{\in C}\), we define the companion sequences
\(G_{j}(q)\) and \(V_{j}(q)\) and a sequence \(W_{j}\) with arbitrary
complex initial values \(w_{0},w_{1}\), by

\begin{equation}
\begin{aligned}
 & G_{0}(q) = 0,\quad G_{1}(q) = 1,\quad V_{0}(q) = 2,\quad V_{1}(q) = 1, \\
 & X_{j + 2} = X_{j + 1} - qX_{j}\quad(X = G,V,W),\quad\quad W_{0} = w_{0},\quad W_{1} = w_{1}.
\end{aligned}
\tag{14}\label{eq:14}
\end{equation}

For fixed \(j\), these quantities are polynomials in \(q\) and \(W_{j}\)
is linear in \(w_{0},w_{1}\). At \(q = - 1\), the companions reduce to
\(F_{j}\) and \(L_{j}\). We consider

\begin{equation}
h = 1 - q,\quad\quad a_{j} = W_{2j} + W_{j} - w_{0}V_{j}(q).
\tag{15}\label{eq:15}
\end{equation}

In particular, \(a_{0} = 0\). The subtraction in (15) is needed when
\(w_{0} \neq 0\); simply replacing \(F_{j}\) by an arbitrary recurrence
sequence in \(F_{j} + F_{2j}\) would not give the stated identity.

\textbf{Theorem 3.} Consider \(q,w_{0},w_{1}\mathbb{\in C}\) and define
\(h,a_{j}\) by (14)--(15). For every even integer \(N \geq 2\):

\begin{equation}
\sum_{j = 0}^{N}\binom{N}{j}h^{N - j}B_{N - j}a_{j} = 0.
\tag{16}\label{eq:16}
\end{equation}

And consequently, for every even integer \(N \geq 4\):

\begin{equation}
N\sum_{j = 1}^{N - 2}\binom{N - 1}{j}h^{N - j}a_{j}\zeta(j + 1 - N) = a_{N} - \frac{Nh}{2}a_{N - 1}.
\tag{17}\label{eq:17}
\end{equation}

Both identities hold at \(q = 1/4\), where the characteristic roots
coincide and at \(q = 1\), where \(h = 0\).

\textbf{Proof.} First we suppose \(q \neq 1/4\) and let \(r,s\) be the
two distinct roots of \(t^{2} - t + q = 0\). Their order is immaterial.
They satisfy \(r + s = 1\), \(rs = q\), and

\begin{equation}
r^{2} + s = 1 - q = h,\quad\quad s^{2} + r = h.
\tag{18}\label{eq:18}
\end{equation}

The recurrence and its initial data give
\(G_{j}(q) = \left( r^{j} - s^{j} \right)/(r - s)\) and
\(V_{j}(q) = r^{j} + s^{j}\). Also,

\begin{equation}
W_{j} = uG_{j}(q) + vV_{j}(q),\quad\quad u = w_{1} - \frac{w_{0}}{2},\quad v = \frac{w_{0}}{2}.
\tag{19}\label{eq:19}
\end{equation}

The equality follows at \(j = 0,1\) and then for all \(j\) by the
recurrence. Substituting (19) into (15), its exponential generating
function becomes

\begin{equation}
\begin{aligned}
 & A(z) = \frac{u}{r - s}\left( e^{r^{2}z} + e^{rz} - e^{s^{2}z} - e^{sz} \right) \\
 & \quad\quad + v\left( e^{r^{2}z} + e^{s^{2}z} - e^{rz} - e^{sz} \right).
\end{aligned}
\tag{20}\label{eq:20}
\end{equation}

Under the operation \(A(z) \mapsto e^{hz}A( - z)\), relation (18)
interchanges the exponents \(r^{2}\) and \(s\) and interchanges
\(s^{2}\) and \(r\). In each bracket of (20), the coefficients at an
interchanged pair have opposite signs. Therefore
\(e^{hz}A( - z) = - A(z)\). Lemma 1 proves (16) for \(q \neq 1/4\). For
fixed \(N,w_{0},w_{1}\), the left-hand side of (16) is a polynomial in
\(q\). It vanishes for every \(q \neq 1/4\), so it is the zero
polynomial and also vanishes at \(q = 1/4\). This argument avoids
division by coincident roots and establishes the repeated-root case
directly at the level of the finite identity.

To derive (17), we isolate the \(j = N\) and \(j = N - 1\) terms of
(16). Since \(a_{0} = 0\), \(B_{0} = 1\) and \(B_{1} = - 1/2\),

\begin{equation}
\sum_{j = 1}^{N - 2}\binom{N}{j}h^{N - j}B_{N - j}a_{j} = \frac{Nh}{2}a_{N - 1} - a_{N}.
\tag{21}\label{eq:21}
\end{equation}

Every Bernoulli index in the remaining sum is at least two, so (3)
applies without its exceptional endpoint. By using the identities

\begin{equation}
B_{N - j} = - (N - j)\zeta(j + 1 - N),\quad\quad(N - j)\binom{N}{j} = N\binom{N - 1}{j},
\tag{22}\label{eq:22}
\end{equation}

substitution of (22) into (21) gives (17) and no division by \(h\)
occurs. At \(h = 0\), equation (16) reads \(a_{N} = 0\), and every
summand in (17) has a positive power of \(h\), so the second identity
remains valid.

\textbf{Corollary 4.} For each integer \(k \geq 2\), the Fibonacci and
Lucas identities are

\begin{equation}
\begin{aligned}
 & 2k\sum_{n = 0}^{2k - 3}\binom{2k - 1}{n + 1}\left( F_{n + 1} + F_{2n + 2} \right)2^{2k - n - 2}\zeta(n + 2 - 2k) \\
 & = \frac{F_{4k} + F_{2k}}{2} - k\left( F_{4k - 2} + F_{2k - 1} \right),
\end{aligned}
\tag{23}\label{eq:23}
\end{equation}

and

\begin{equation}
\begin{aligned}
 & 2k\sum_{n = 0}^{2k - 3}\binom{2k - 1}{n + 1}\left( L_{2n + 2} - L_{n + 1} \right)2^{2k - n - 2}\zeta(n + 2 - 2k) \\
 & = \frac{L_{4k} - L_{2k}}{2} - k\left( L_{4k - 2} - L_{2k - 1} \right).
\end{aligned}
\tag{24}\label{eq:24}
\end{equation}

\textbf{Proof.} In (17), we put \(q = - 1\), \(h = 2\), \(N = 2k\) and
divide both sides by two. Choosing
\(\left( w_{0},w_{1} \right) = (0,1)\) gives \(a_{j} = F_{j} + F_{2j}\);
choosing \(\left( w_{0},w_{1} \right) = (2,1)\) gives
\(a_{j} = L_{2j} - L_{j}\). Finally \(j = n + 1\).

For an arbitrary solution of the Fibonacci recurrence, the same proof
uses \(a_{j} = W_{2j} + W_{j} - W_{0}L_{j}\). Equally, it combines (23)
and (24) with coefficients \(w_{1} - w_{0}/2\) and \(w_{0}/2\).

This relationship with previous convolution formulas is explicit which
for \(h \neq 0\) and an exponential weight \(\lambda^{j}\), equation (2)
gives

\begin{equation}
T_{N}^{(h)}\left( \left( \lambda^{j} \right)_{j \geq 0} \right) = h^{N}B_{N}(\lambda/h).
\tag{25}\label{eq:25}
\end{equation}

Applying (25) to the four exponentials in (20) and then reflecting their
paired arguments gives (16). For \(q = - 1\), the Fibonacci and Lucas
evaluations used in this calculation are the specializations \(x = 0\),
\(z = 1/2\), \(m = 0\). Thus Corollary 4 lies within the earlier
Bernoulli-polynomial frame. Theorem 3 records the recurrence parameter
and the exceptional cases in a single polynomial identity. Two limits
clarify the scope. At \(q = 1/4\), the repeated-root formulas are
\(G_{j}(1/4) = j2^{1 - j}\) and \(V_{j}(1/4) = 2^{1 - j}\), including
the prescribed values at \(j = 0\). Here \(h = 3/4\) and Theorem 3
remains a nontrivial finite identity. At \(q = 1\), the scaling
parameter is zero and the theorem reduces to the even-index relation
\(a_{N} = 0\). The identity is therefore valid there, but its zeta sum
no longer contains nonzero terms. This construction uses the condition
\(r + s = 1\). For a recurrence with distinct roots \(r,s\) of arbitrary
sum, the two desired pairings would require \(r^{2} + s = s^{2} + r\).
Their difference is \((r - s)(r + s - 1)\), so the same pairing forces
\(r + s = 1\).

\hypertarget{absolute-sums}{%
\section*{\texorpdfstring{3 Absolute sums
}{3 Absolute sums }}\label{absolute-sums}}

The signed identities in Section 2 depend on reflection. Their
absolute-term analysis admits a wider class of weights. Throughout this
section, we let \(\left( a_{j} \right)_{j \geq 0}\) be a real sequence
such that \(a_{0} = 0\), \(a_{j} \geq 0\) for \(j \geq 1\),
\(a_{2} > 0\), and \(a_{j} \leq C\rho^{j}\) for some constants
\(C,\rho > 0\). We define the even exponential generating function and
its coefficients at \(\pi\) by

\begin{equation}
E(t) = \sum_{\substack{j \geq 2 \\ j\ even}}^{}a_{j}\frac{t^{j}}{j!},\quad\quad w_{j} = \frac{a_{j}\pi^{j}}{j!}\quad\left( j \geq 2\text{ even} \right),\quad\quad D = E(\pi).
\tag{26}\label{eq:26}
\end{equation}

The growth assumption makes the series in (26) entire and \(E(t) > 0\)
for real \(t > 0\). In particular, \(D > 0\). For an even integer
\(N \geq 4\), by using the original zeta normalization and setting

\begin{equation}
\begin{aligned}
 & t_{N,j} = N\binom{N - 1}{j}a_{j}2^{N - j - 1}\zeta(j + 1 - N),\quad 1 \leq j \leq N - 2, \\
 & M_{N} = \sum_{j = 1}^{N - 2}{|t_{N,j}|},\quad\quad S_{N} = \frac{\pi^{N}M_{N}}{N!},
\end{aligned}
\tag{27}\label{eq:27}
\end{equation}

odd indices in (27) give zero terms since their zeta arguments are
negative even integers. At an even index, we put \(m = N - j \geq 2\).
The classical formula \(|B_{m}| = 2m!\zeta(m)/(2\pi)^{m}\), together
with (3), gives the exact reduction

\begin{equation}
|t_{N,j}| = \frac{N!}{\pi^{N}}w_{j}\zeta(N - j),\quad\quad S_{N} = \sum_{\substack{2 \leq j \leq N - 2 \\ j\ even}}^{}w_{j}\zeta(N - j).
\tag{28}\label{eq:28}
\end{equation}

The first equality in (28) is for even \(j\). It replaces alternating
negative-integer zeta values by a finite sum with positive coefficients
and positive-integer zeta values. The ordinary convergent Dirichlet
series is therefore available throughout the following proof.

\hypertarget{a-finite-exponential-expansion}{%
\subsection*{3.1 A finite exponential
expansion}\label{a-finite-exponential-expansion}}

\textbf{Theorem 5.} Under the assumptions of this section, for every
integer \(K \geq 1\) and even integers \(N \geq 4\),

\begin{equation}
S_{N} = \sum_{m = 1}^{K}\frac{E(m\pi)}{m^{N}} + R_{N,K},\quad\quad|R_{N,K}| \leq \frac{(K + 2)E\left( (K + 1)\pi \right)}{(K + 1)^{N}}.
\tag{29}\label{eq:29}
\end{equation}

For each fixed integer \(K \geq 0\), we define a scaled remainder with
an empty sum when \(K = 0\), by

\begin{equation}
H_{N,K} = \frac{(K + 1)^{N}}{E\left( (K + 1)\pi \right)}\left( S_{N} - \sum_{m = 1}^{K}\frac{E(m\pi)}{m^{N}} \right).
\tag{30}\label{eq:30}
\end{equation}

For the remainders in (30), \(H_{N,K} \rightarrow 1\) as \(N\) tends to
infinity through even integers and the following bound holds for every
even \(N \geq 4\):

\begin{equation}
|H_{N,K} - 1| \leq \frac{(K + 3)E\left( (K + 2)\pi \right)}{E\left( (K + 1)\pi \right)}\left( \frac{K + 1}{K + 2} \right)^{N}.
\tag{31}\label{eq:31}
\end{equation}

\textbf{Proof.} First, we separate the first \(K\) terms of the
Dirichlet series for each zeta value in (28). Completing each resulting
finite sum in \(j\) to the entire function \(E\) adds a tail that must
then be subtracted. This gives us \(R_{N,K} = U_{N,K} - V_{N,K}\), where

\begin{equation}
\begin{aligned}
 & U_{N,K} = \sum_{\substack{2 \leq j \leq N - 2 \\ j\ even}}^{}w_{j}\sum_{m = K + 1}^{\infty}m^{- (N - j)}, \\
 & V_{N,K} = \sum_{m = 1}^{K}m^{- N}\sum_{\substack{j \geq N \\ j\ even}}^{}a_{j}\frac{(m\pi)^{j}}{j!}.
\end{aligned}
\tag{32}\label{eq:32}
\end{equation}

Both quantities are nonnegative and the outer sum defining \(U_{N,K}\)
is finite and each inner series has exponent at least two; the sums in
\(V_{N,K}\) converge by the exponential growth assumption. These facts
justify the rearrangement leading to (32). For every real \(s \geq 2\),
monotonicity of \(x^{- s}\) yields

\begin{equation}
\sum_{m = K + 1}^{\infty}m^{- s} \leq (K + 1)^{- s} + \int_{K + 1}^{\infty}x^{- s}\, dx \leq (K + 2)(K + 1)^{- s}.
\tag{33}\label{eq:33}
\end{equation}

Applying (33) with \(s = N - j\) and completing the positive sum in
\(j\) gives

\begin{equation}
0 \leq U_{N,K} \leq \frac{(K + 2)E\left( (K + 1)\pi \right)}{(K + 1)^{N}}.
\tag{34}\label{eq:34}
\end{equation}

Considering \(j \geq N\) and \(1 \leq m \leq K\), the exponent \(j - N\)
is nonnegative, so \(m^{j - N} \leq (K + 1)^{j - N}\). Substituting this
inequality into the second line of (32) gives

\begin{equation}
0 \leq V_{N,K} \leq \frac{KE\left( (K + 1)\pi \right)}{(K + 1)^{N}}.
\tag{35}\label{eq:35}
\end{equation}

Hence \(|U - V| \leq max(U,V)\) for nonnegative \(U,V\), equations (34)
and (35) prove Eq.(29). We apply that result with \(K + 1\) in place of
\(K\) and multiply the remainder by
\((K + 1)^{N}/E\left( (K + 1)\pi \right)\). The first unremoved term
becomes one, giving us (31). Its right-hand side tends to zero for every
fixed \(K\) which proves the limit.

For \(K = 1\), equation (29) gives \(|S_{N} - D| \leq 3E(2\pi)2^{- N}\).
Higher fixed values of \(K\) identify successive contributions with
exponential scales \(2^{- N},3^{- N},\ldots\). The constants depend on
\(K\) and the weight sequence; no uniform assertion is made when \(K\)
grows with \(N\). In particular, Eq. (29) is a finite asymptotic
expansion and does not assert convergence of an infinite series in
\(m\).

\hypertarget{convergence-of-the-term-magnitudes}{%
\subsection*{3.2 Convergence of the term
magnitudes}\label{convergence-of-the-term-magnitudes}}

We normalize the nonzero term magnitudes to obtain a probability
distribution on the positive even integers:

\begin{equation}
p_{N,j} = \frac{|t_{N,j}|}{M_{N}}\quad(2 \leq j \leq N - 2),\quad\quad p_{j} = \frac{w_{j}}{D}\quad\left( j \geq 2\text{ even} \right).
\tag{36}\label{eq:36}
\end{equation}

We extend \(p_{N,j}\) by zero for every even \(j \geq N\). Both
distributions in (36) then sum to one. They describe deterministic
shares of the absolute-term sum. For distributions on this countable
set, the total variation distance is one half of the sum of the absolute
coordinate differences.

\Needspace{6\baselineskip}
\textbf{Theorem 6.} Under the assumptions, for every even integer
\(N \geq 4\),

\begin{equation}
\frac{1}{2}\sum_{\substack{j \geq 2 \\ j\ even}}^{}{|p_{N,j} - p_{j}|} \leq min\{ 1,\frac{3E(2\pi)}{D}\, 2^{- N}\}.
\tag{37}\label{eq:37}
\end{equation}

Thus, \(p_{N,j} \rightarrow p_{j}\) uniformly in the even index \(j\).

\textbf{Proof.} We consider \(u_{N,j} = w_{j}\zeta(N - j)\) for
\(2 \leq j \leq N - 2\) and \(u_{N,j} = 0\) for \(j \geq N\), for even
indices only. Equation (28) gives \(\sum_{j}^{}u_{N,j} = S_{N}\). The
excess on the retained indices and the omitted tail are

\begin{equation}
I_{N} = \sum_{\substack{2 \leq j \leq N - 2 \\ j\ even}}^{}w_{j}\lbrack\zeta(N - j) - 1\rbrack,\quad\quad J_{N} = \sum_{\substack{j \geq N \\ j\ even}}^{}w_{j}.
\tag{38}\label{eq:38}
\end{equation}

The quantities in (38) satisfy \(S_{N} - D = I_{N} - J_{N}\) and
\(\sum_{j}^{}{|u_{N,j} - w_{j}|} = I_{N} + J_{N}\). The bound for
\(U_{N,1}\) in (34) applies to \(I_{N}\), while \(2^{j} \geq 2^{N}\) on
the tail gives

\begin{equation}
I_{N} \leq 3E(2\pi)2^{- N},\quad\quad J_{N} \leq E(2\pi)2^{- N}.
\tag{39}\label{eq:39}
\end{equation}

To control the normalization, we write
\[u_{N,j}/S_{N} - w_{j}/D = \left( u_{N,j} - w_{j} \right)/D + \left( D - S_{N} \right)u_{N,j}/\left( DS_{N} \right).\]
Summing absolute values and using \(\sum_{j}^{}u_{N,j} = S_{N}\) yields

\begin{equation}
\frac{1}{2}\sum_{j}^{}{|\frac{u_{N,j}}{S_{N}} - \frac{w_{j}}{D}|} \leq \frac{I_{N} + J_{N} + |I_{N} - J_{N}|}{2D} = \frac{\max\left( I_{N},J_{N} \right)}{D}.
\tag{40}\label{eq:40}
\end{equation}

Equations (39) and (40) prove the nontrivial bound in (37). The distance
between probability distributions is also at most one. A single
coordinate difference is bounded by total variation, since the positive
and negative coordinate differences have equal total mass. This proves
uniform convergence. The limiting distribution actually depends on the
weights through \(a_{j}\pi^{j}/j!\). The factorial in the denominator
controls its large-index tail. A strict largest coefficient of this
sequence is therefore a natural candidate for the largest finite-sum
term.

\hypertarget{fibonacci-and-lucas-cancellation-profiles}{%
\section*{4 Fibonacci and Lucas cancellation
profiles}\label{fibonacci-and-lucas-cancellation-profiles}}

\hypertarget{closed-values-and-conditioning}{%
\subsection*{4.1 Closed values and
conditioning}\label{closed-values-and-conditioning}}

Considering the weights in Corollary 4 and distinguish them by a
superscript \(X \in \{ F,L\}\):

\begin{equation}
a_{j}^{F} = F_{j} + F_{2j},\quad\quad a_{j}^{L} = L_{2j} - L_{j},\quad\quad C_{N}^{X} = \frac{a_{N}^{X} - Na_{N - 1}^{X}}{2}.
\tag{41}\label{eq:41}
\end{equation}

We use the label \(X\) on
\(t_{N,j},M_{N},S_{N},E,D,w_{j},p_{N,j},p_{j},H_{N,K}\) when those
objects are formed from \(a_{j}^{X}\). The identities (23) and (24) say
that \(\sum_{j}^{}t_{N,j}^{X} = C_{N}^{X}\) with \(N = 2k\). For
\(k = 2,3,4,5\), the Fibonacci values are \(- 8, - 104, - 1056, - 9680\)
and the Lucas values are \(- 8, - 184, - 2176, - 21008\).

All even weights are positive and \(a_{2}^{F} = a_{2}^{L} = 4\). By
putting \(r = \alpha^{2}\), the estimates \(F_{j} \leq \alpha^{j}\) and
\(L_{2j} = r^{j} + \beta^{2j} \leq 2r^{j}\) imply
\(0 \leq a_{j}^{X} \leq 2r^{j}\) for \(j \geq 1\). Here
\(a_{1}^{L} = 2\), and positivity for larger indices follows from
monotonicity of the Lucas sequence from index one. Thus both sequences
satisfy the assumptions of Section 3. Taking even parts in (8) evaluates
the required entire functions are

\begin{equation}
\begin{aligned}
E_{F}(t) & = \frac{\cosh(\alpha t) + cosh\left( \alpha^{2}t \right) - cosh(\beta t) - cosh\left( \beta^{2}t \right)}{\sqrt{5}}, \\
E_{L}(t) & = cosh\left( \alpha^{2}t \right) + cosh\left( \beta^{2}t \right) - cosh(\alpha t) - cosh(\beta t).
\end{aligned}
\tag{42}\label{eq:42}
\end{equation}

Here we have \(cosht = \left( e^{t} + e^{- t} \right)/2\). Equation (42)
gives \(D_{F} = E_{F}(\pi) \approx 868.247964\) and
\(D_{L} = E_{L}(\pi) \approx 1783.786571\). The signed totals never
vanish at the degrees under consideration. Direct use of the recurrence
gives, for even \(N \geq 4\),

\begin{equation}
\begin{aligned}
3a_{N - 1}^{F} - a_{N}^{F} & = F_{2N - 4} + F_{N - 1} + F_{N - 3} > 0, \\
3a_{N - 1}^{L} - a_{N}^{L} & = L_{2N - 4} - L_{N - 1} - L_{N - 3} > 0.
\end{aligned}
\tag{43}\label{eq:43}
\end{equation}

For the second line, \(2N - 4 \geq N\) and
\(L_{N} = L_{N - 1} + L_{N - 2} > L_{N - 1} + L_{N - 3}\). Since
\(N > 3\) and \(a_{N - 1}^{X} > 0\), (43) implies \(C_{N}^{X} < 0\). The
component wise summation condition number is therefore well defined:

\begin{equation}
\kappa_{N}^{X} = \frac{M_{N}^{X}}{|C_{N}^{X}|}.
\tag{44}\label{eq:44}
\end{equation}

If the individual terms are perturbed by
\(|\delta t_{N,j}^{X}| \leq \varepsilon|t_{N,j}^{X}|\), their sum
changes by at most \(\varepsilon M_{N}^{X}\), so the relative change is
bounded by \(\varepsilon\kappa_{N}^{X}\). Equality is attained by giving
every perturbation the same sign and its allowed magnitude. Thus (44)
measures sensitivity to termwise perturbations; the actual rounding
error of a particular algorithm is a separate quantity.

\hypertarget{a-fixed-largest-index}{%
\subsection*{4.2 A fixed largest index}\label{a-fixed-largest-index}}

\textbf{Proposition 7.} For \(X = F\) and \(X = L\), the limiting
distribution \(p_{j}^{X}\) has its unique maximum at \(j = 8\). For
every even integer \(N \geq 26\), the unique largest term magnitude
among \(|t_{N,j}^{X}|\), \(1 \leq j \leq N - 2\), also occurs at
\(j = 8\). As \(N\) tends to infinity through even integers,

\begin{equation}
\max_{1 \leq j \leq N - 2}|t_{N,j}^{X}| \sim \frac{a_{8}^{X}\pi^{8}}{8!}\frac{N!}{\pi^{N}},\quad\quad p_{N,8}^{X} \rightarrow \frac{a_{8}^{X}\pi^{8}}{8!D_{X}}.
\tag{45}\label{eq:45}
\end{equation}

The limiting shares in (45) are approximately \(0.273209\) for Fibonacci
weights and \(0.284964\) for Lucas weights. The threshold 26 is
sufficient and no minimal-threshold assertion is made.

\textbf{Proof.} We first find the largest two limiting weights, then use
Theorem 6 to preserve their gap at finite degree. The recurrence gives

\begin{equation}
\begin{aligned}
\left( a_{2}^{F},a_{4}^{F},a_{6}^{F},a_{8}^{F},a_{10}^{F} \right) & = (4,24,152,1008,6820), \\
\left( a_{2}^{L},a_{4}^{L},a_{6}^{L},a_{8}^{L},a_{10}^{L} \right) & = (4,40,304,2160,15004).
\end{aligned}
\tag{46}\label{eq:46}
\end{equation}

For even \(j\), the ratio of consecutive limiting weights is

\begin{equation}
\frac{w_{j + 2}^{X}}{w_{j}^{X}} = \frac{a_{j + 2}^{X}}{a_{j}^{X}}\frac{\pi^{2}}{(j + 1)(j + 2)}.
\tag{47}\label{eq:47}
\end{equation}

Using \(9 < \pi^{2} < 10\) in (46)--(47) gives
\(w_{2}^{X} < w_{4}^{X} < w_{6}^{X} < w_{8}^{X}\). It also gives
\(w_{10}^{X} < w_{6}^{X}\), since \(100a_{10}^{X} < 5040a_{6}^{X}\) for
each of the two listed sequences.

We next show \(a_{j + 2}^{X} < 8a_{j}^{X}\) for every even \(j \geq 8\).
For Fibonacci numbers, \(F_{n + 4} = 5F_{n} + 3F_{n - 1} < 8F_{n}\) for
\(n \geq 3\), and \(F_{j + 2} \leq 3F_{j}\), so the assertion follows
after adding the terms at indices \(j + 2\) and \(2j + 4\). For Lucas
numbers, \(L_{2j + 4} = 7L_{2j} - L_{2j - 4} < 7L_{2j}\). Also
\(L_{2j} = L_{j}^{2} - 2 > 8L_{j}\) for even \(j \geq 8\), because
\(L_{j} \geq L_{8} = 47\). Hence \(a_{j}^{L} > 7L_{2j}/8\) and
\(a_{j + 2}^{L} < L_{2j + 4} < 7L_{2j} < 8a_{j}^{L}\). The identities
used here follow by substitution in (8), or directly from the
recurrence. Eq. (47) now gives
\(w_{j + 2}^{X}/w_{j}^{X} < 80/\lbrack(j + 1)(j + 2)\rbrack < 1\) for
every even \(j \geq 8\). Therefore \(w_{8}^{X}\) is the unique largest
weight and \(w_{6}^{X}\) is the second largest. Their gap has the
explicit lower bound

\begin{equation}
w_{8}^{X} - w_{6}^{X} = \frac{\pi^{6}}{8!}\left( a_{8}^{X}\pi^{2} - 56a_{6}^{X} \right) > \frac{729}{40320}\left( 9a_{8}^{X} - 56a_{6}^{X} \right) \geq \frac{81}{8}.
\tag{48}\label{eq:48}
\end{equation}

To transfer this gap to \(p_{N,j}^{X}\), Theorem 6 bounds each
coordinate error by \(3E_{X}(2\pi)2^{- N}/D_{X}\). For any even
\(j \neq 8\),

\begin{equation}
p_{N,8}^{X} - p_{N,j}^{X} \geq \frac{w_{8}^{X} - w_{6}^{X} - 6E_{X}(2\pi)2^{- N}}{D_{X}}.
\tag{49}\label{eq:49}
\end{equation}

It remains to make the numerator positive uniformly in \(N\). The
exponential weight bound gives
\(E_{X}(2\pi) \leq 2\lbrack cosh(2\pi r) - 1\rbrack < 2e^{2\pi r}\).
Since \(r < 8/3\), \(\pi < 22/7\) and \(e < 11/4\), we have
\(2\pi r < 17\). So, for \(N \geq 26\),

\begin{equation}
6E_{X}(2\pi)2^{- N} < \frac{12(11/4)^{17}}{2^{26}} = \frac{15163410854978}{28823037615171} < \frac{81}{8}.
\tag{50}\label{eq:50}
\end{equation}

Combining Eqs. (48)--(50) proves the strict inequality in (49).
Odd-indexed terms are zero, so \(j = 8\) is also the unique largest
index among all summands. Finally, (28) at \(j = 8\) and
\(\zeta(N - 8) \rightarrow 1\) prove the first asymptotic in (45) and
Theorem 6 gives the limiting share. In the original indexing of
Corollary 4, the largest term is at \(n = j - 1 = 7\). A threshold is
necessary for this argument: at \(N = 12\), direct exact comparison
gives the largest index \(j = 10\) for both sequences. Proposition 7 is
a statement about every even \(N \geq 26\), rather than an extrapolation
from finitely many computed maxima.

\hypertarget{growth-of-the-cancellation}{%
\subsection*{4.3 Growth of the
cancellation}\label{growth-of-the-cancellation}}

\textbf{Corollary 8.} We let \(c_{F} = 1/\sqrt{5}\), \(c_{L} = 1\) and
\(r = \alpha^{2}\). As \(N\) tends to infinity through even integers,

\begin{equation}
M_{N}^{X} \sim \frac{N!}{\pi^{N}}D_{X},\quad\quad\kappa_{N}^{X} \sim \frac{2rD_{X}}{c_{X}}\frac{(N - 1)!}{(\pi r)^{N}}.
\tag{51}\label{eq:51}
\end{equation}

A more accurate denominator retains the finite shift \(r\):

\begin{equation}
\kappa_{N}^{X} = \frac{2rD_{X}}{c_{X}}\frac{N!}{(N - r)(\pi r)^{N}}\lbrack 1 + O\left( \alpha^{- N} \right)\rbrack.
\tag{52}\label{eq:52}
\end{equation}

The constants implicit in the last estimate depend only on \(X\).

\textbf{Proof.} The case \(K = 1\) of (29) gives
\(M_{N}^{X} = \left( N!/\pi^{N} \right)\lbrack D_{X} + O\left( 2^{- N} \right)\rbrack\).
Binet's formulas give
\(a_{n}^{X} = c_{X}r^{n}\lbrack 1 + O\left( \alpha^{- n} \right)\rbrack\)
which the largest subleading exponential has base \(\alpha\) and
\(\alpha/r = \alpha^{- 1}\). Substitution in (41), using
\(C_{N}^{X} < 0\), yields

\begin{equation}
|C_{N}^{X}| = \frac{c_{X}}{2}r^{N - 1}(N - r)\lbrack 1 + O\left( \alpha^{- N} \right)\rbrack.
\tag{53}\label{eq:53}
\end{equation}

Indeed, the two relative remainder contributions are bounded by constant
multiples of \(\alpha^{- N}r/(N - r)\) and
\(N\alpha^{- (N - 1)}/(N - r)\), respectively, and both have order
\(\alpha^{- N}\). Dividing the absolute-sum estimate by (53). Since
\(2^{- N} = O\left( \alpha^{- N} \right)\), this gives us (52). Finally,
\(N/(N - r) \rightarrow 1\), which proves (51). Here
\(u_{N} \sim v_{N}\) means \(u_{N}/v_{N} \rightarrow 1\), and
\(O\left( \alpha^{- N} \right)\) denotes a quantity bounded in absolute
value by a constant times \(\alpha^{- N}\) for all sufficiently large
even \(N\).

The absolute-term distribution therefore stabilizes while the condition
number grows factorially after division by a fixed exponential factor.
These are compatible facts: a fixed collection of early indices carries
most of an absolute-term sum whose scale is \(N!/\pi^{N}\), whereas its
signed total has only recurrence-scale growth. The finite-order
restriction in Theorem 5 is substantive for these weights. Eq. (42)
gives \(E_{X}(m\pi) \sim \left( c_{X}/2 \right)e^{\pi rm}\) as
\(m \rightarrow \infty\). Accordingly, for any fixed \(N\), the terms
\(E_{X}(m\pi)/m^{N}\) do not tend to zero. The infinite series suggested
by extending (29) to all \(m\) diverges. Each fixed truncation remains
valid with the stated remainder bound as \(N\) increases.

In \textbf{Figure 1}, the intermediate cancellation at \(N = 48\)
(\(k = 24\)) is examined, by defining
\(P_{n}^{X} = \sum_{j = 1}^{n + 1}t_{48,j}^{X}\) for
\(0 \leq n \leq 45\). The normalized partial sums have excursions of
order \(10^{18}\), although their exact final value is one. That
endpoint is below the graph's visual resolution on this scale; its value
comes from exact arithmetic.

\begin{figure}[!htbp]
\centering
% Original image alternative text: Figure 1. Normalized partial sums P_n^X/C_{48}^X for the Fibonacci and Lucas identities at N=48, with markers calculated from exact rational sums and connecting lines used as visual guides.
\includegraphics[width=0.92\linewidth]{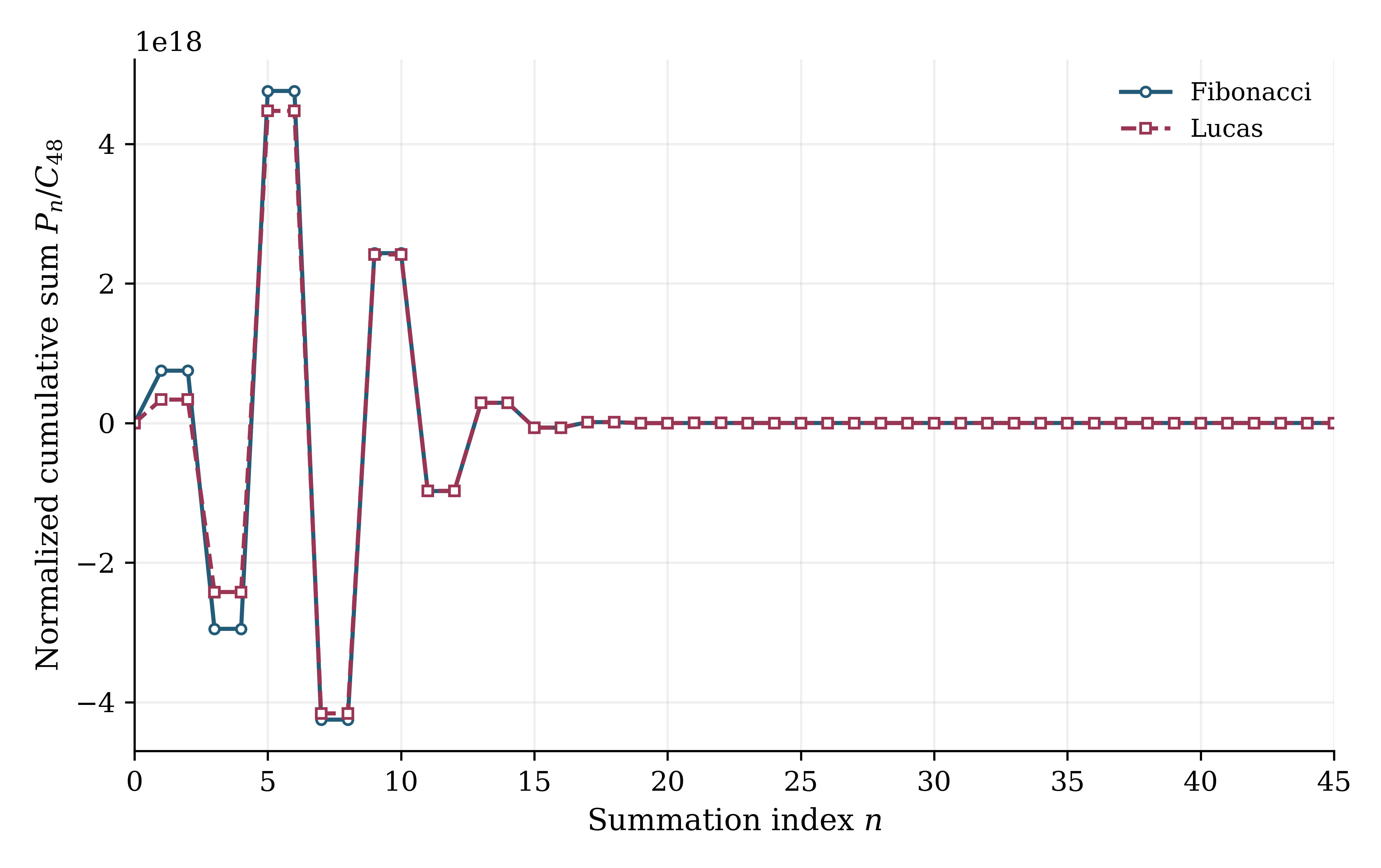}
\caption{Normalized partial sums \(P_{n}^{X}/C_{48}^{X}\) for
the Fibonacci and Lucas identities at \(N = 48\).}
\label{fig:1}
\end{figure}

The way that how the absolute-term sum is distributed across its indices
is shown in \textbf{Figure 2}. The finite-degree markers at \(N = 48\)
visually coincide with the limiting masses in (36). The total variation
distances are near \(9.74 \times 10^{- 12}\) for Fibonacci weights and
\(1.04 \times 10^{- 11}\) for Lucas weights. Only even indices are shown
since every odd-indexed summand is zero.

\begin{figure}[!htbp]
\centering
% Original image alternative text: Figure 2. Shares p_{48,j}^X of the absolute-term sum at even indices 2\leq j\leq30, shown as exact-rational computed markers against the limiting masses p_j^X joined by visual guides, with circles and a solid line for Fibonacci weights and squares and a dashed line for Lucas weights.
\includegraphics[width=0.92\linewidth]{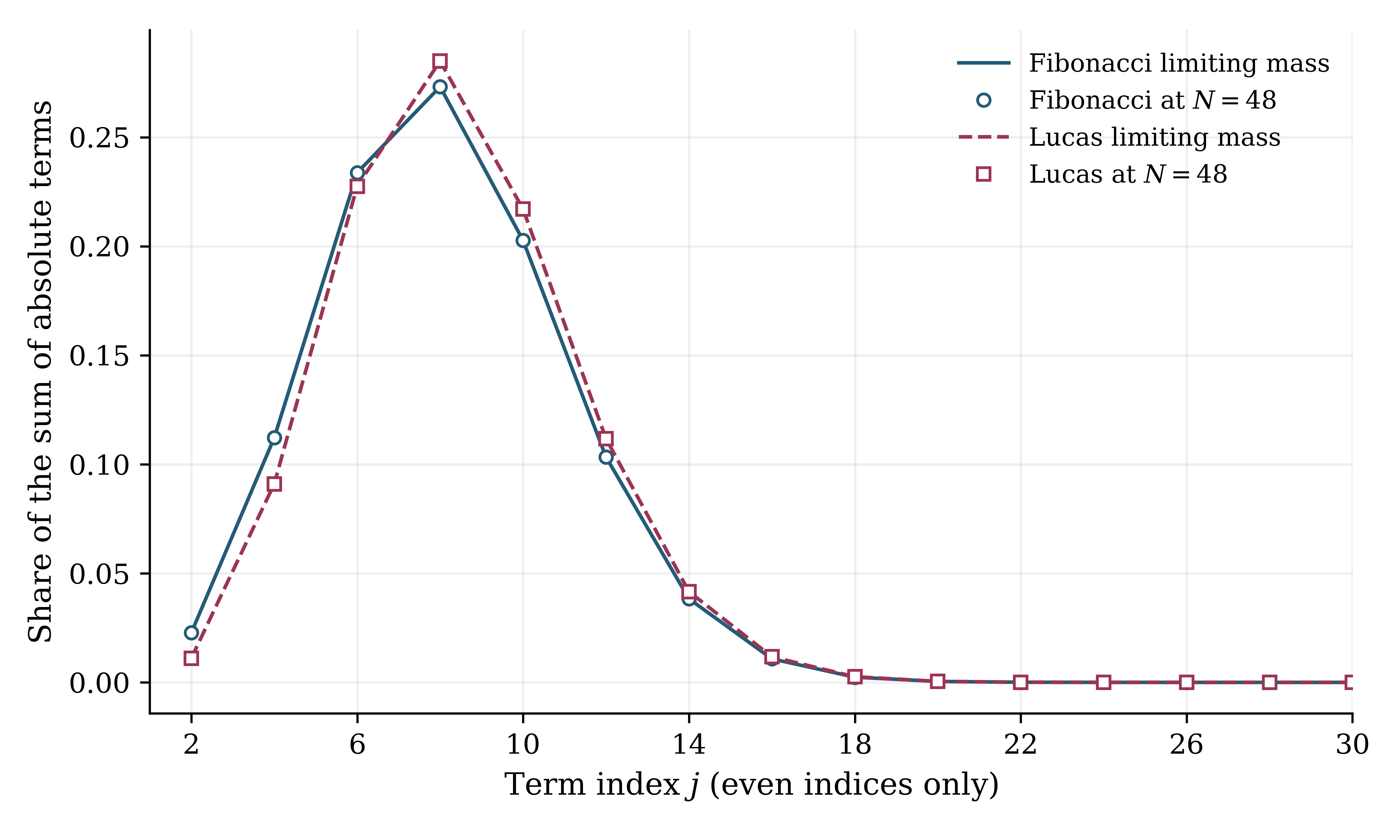}
\caption{Shares \(p_{48,j}^{X}\) of the absolute-term sum at
even indices \(2 \leq j \leq 30\).}
\label{fig:2}
\end{figure}

\FloatBarrier
\hypertarget{discussion}{%
\section*{\texorpdfstring{5 Discussion
}{5 Discussion }}\label{discussion}}

Our main contribution in this work, is a quantitative description of how
the magnitudes of the cancelling terms are distributed across the two
finite zeta sums. Although the number of terms grows with the degree, we
prove that the largest term in absolute value remains at index eight for
every even degree of at least twenty-six. An explicit limiting
distribution, a total variation bound and a strict gap between its two
largest coefficients make this a rigorous finite-degree result. Our
analysis develops the absolute-term expansion and limiting magnitude
profile under a broader assumption: nonnegative weights with at most
exponential growth. These results explain how exact, compact expressions
can coexist with severe sensitivity in direct summation, while
recurrence-based evaluation avoids forming the large cancelling terms.
The expansion is valid at each fixed truncation order, although its
infinite continuation diverges for the Fibonacci and Lucas weights. The
degree threshold is sufficient rather than optimal and neither
limitation weakens the stated theorem.

\medskip

Department of Energy, Politecnico di Torino, Italy\\
\href{mailto:payam.danesh@polito.it}{\nolinkurl{payam.danesh@polito.it}}

\end{document}